\documentclass[graybox]{svmult}

\usepackage{type1cm}        
\usepackage{makeidx}         
\usepackage{graphicx}        
\usepackage{multicol}        
\usepackage[bottom]{footmisc}

\usepackage{newtxtext}       %
\usepackage[varvw]{newtxmath}       

\makeindex             

\usepackage{dsfont}
\usepackage{pgffor}
\usepackage{booktabs}
\usepackage{subcaption}
\usepackage{algpseudocode}

\renewcommand{\I}{\mathcal{I}}

\renewcommand{\S}{\mathcal{S}}

\newcommand{\polyspace}{\mathds{P}}

\newcommand{\Real}{\mathbb{R}}

\newcommand{\domain}{\Omega}
\newcommand{\hadomain}{\tilde{\Omega}}
\newcommand{\Tau}{\mathcal{T}}
\newcommand{\mesh}{\Tau_h}
\newcommand{\cells}{\mathcal{K}}

\newcommand{\cell}{K}

\newcommand{\point}{P}

\newcommand{\cellCentre}{\vec{k}}

\newcommand{\pointvec}{\vec{p}}

\newcommand{\triangleCentre}{\vec{t}}

\newcommand{\mat}[1]{\expandafter\MakeUppercase\expandafter{\text{#1}}}
\newcommand{\pdd}[3]{\frac{\partial^{#2} #1}{\partial #3^{#2}}}
\newcommand{\norm}[1]{\left\lVert#1\right\rVert}
\newcommand{\abs}[1]{\left\vert#1\right\vert}

\newcommand{\argmin}[1]{\underset{#1}{\text{argmin}}}

\newcommand{\snGrad}[1][]{\nabla^\perp_{h
    \if#1\else,#1\fi
}}
\newcommand{\faceGrad}[1][]{\nabla_{F
    \if#1\else_{#1}\fi
}}
\newcommand{\corrSnGrad}[1][]{\widetilde{\nabla}^\perp_{h
    \if#1\else_{#1}\fi
}}
\newcommand{\ggGrad}[1][]{\nabla_{h
    \if#1\else,#1\fi
}}
\newcommand{\faceInterp}[1][]{\I_{h
    \if#1\else,#1\fi
}}

\begin{document}

\title*{An Immersed Boundary Method for Polymeric Continuous Mixing}
\author{Giorgio Negrini, Nicola Parolini and Marco Verani}
\institute{Giorgio Negrini \at Pirelli Tyre S.p.a., Viale P. e A. Pirelli, 25, 20126 Milano, \email{giorgio.negrini@pirelli.com}
\and Nicola Parolini \at MOX, Dipartimento di Matematica, Politecnico di Milano, P. Leonardo da Vinci, 32, 20133 Milano \email{nicola.parolini@polimi.it}
\and Marco Verani \at MOX, Dipartimento di Matematica, Politecnico di Milano, P. Leonardo da Vinci, 32, 20133 Milano \email{marco.verani@polimi.it}
}
%
%
\maketitle

\abstract{
We introduce a new implementation of the Immersed Boundary method in the finite-volume library OpenFOAM. The implementation is tailored to the simulation of temperature-dependent non-Newtonian polymeric flows in complex moving geometries, such as those characterizing the most popular polymeric mixing technologies. 
}

\section{Introduction}
Mixing is a fundamental task in many industrial processes and consists in manipulating a heterogeneous physical system, with the intent to make it more homogeneous. In the polymer processing industry, mixing tasks can be performed with different technologies, that can be classified as either \textit{continuous} or \textit{batch} mixers. The latter have to be cyclically emptied and reloaded, while continuous mixers, including single-screw, twin-screw and planetary roller extruders, develop a steady and constant flow of material pushed by one or more screw through a cylindrical barrel. In the past decades, many efforts based both on experimental analysis and numerical simulations have been devoted to their optimization, in particular of continuous mixing devices that proved to be the most effective technology.

Single-screw extruders (SSE) and twin-screw extruders (TSE) are widely used in mixing processes and they have been extensively studied in the literature
\cite{JI2020}. More recently, a new class of extruders, the so-called planetary roller extruders (PRE), characterized by more complex geometries and kinematics, has emerged has a competitive alternative for specific mixing tasks. The PRE is a multi-screw extruder composed by a central spindle (sun), a barrel (ring) and a variable number of smaller spindles (planets) between them. The rotation of the sun drives the one of the planets thanks to their gear-like shapes. The flow in a PRE is driven by both drag forces and pressure gradients because the spindles are characterized by a helical gearing to allow the transport of the fluid along its axis. Moreover, the gearing increases significantly the contact surface between the fluid and the spindles with respect to single- and twin- screw extruders, which allows the drag to noticeably drive the flow. A planetary roller extruder is the optimal choice in many demanding continuous mixing tasks since it is able to efficiently masticate, mix, homogenize, disperse and de-gas highly viscous substrates. Very few numerical explorations have been done on this type of device. Up to the authors knowledge, the only work presenting comprehensive numerical simulations of PRE is \cite{WINCK2021}, where a sector of a PRE is considered and all the wall boundaries (sun, ring and spindles) are approximated using the mesh superposition technique \cite{MST} implemented in ANSYS Polyflow.

In this paper, we propose a new implementation of the Immersed Boundary (IB) method in the open-source finite volume library OpenFOAM, developed with the objective of simulating continuous mixing processes. The choice of resorting to a non-conforming approach such as IB is driven by the impossibility of treating the PRE configuration with body fitted approaches. Indeed, body-fitted methods that have been successfully proposed to deal with both single-screw \cite{SIRJALA2000, CHANG1995} and twin-screw \cite{HELMIG2019, HINZ2020} cannot be adapted to the complex kinematics of PRE.  

Different non-conforming discretization approaches have been proposed in the literature in the past decades. With no claim of being exhaustive we mention the Immersed Boundary method \cite{PESKIN2002,MITTAL2005,BOFFI2003}, the Fictitious Domain \cite{GLOWINSKI1994283,massing2014stabilized} method, the Diffuse Interface method \cite{SCHLOTTBOM2014,SCHLOTTBOM2016,NEGRINI2021}, the eXtended Finite Element method \cite{SCHOTT2014233}, the Mesh Superposition Technique \cite{MST}, the Cut Finite Element method \cite{BURMAN2010,BURMAN2012} that have been presented in different flavours and in multiple discretization frameworks (finite differences, finite volumes, finite elements).  
Here, we revisit the IBM proposed in \cite{JASAKIBM2014} in order to improve its performance in terms of accuracy, in particular when applied to peculiar geometries such as those characterizing continuous mixing devices. 
Moreover, the scalability performances has been improved allowing the solution of large cases of industrial interest.

The paper is organized as follows: in Section \ref{sec:model} the non-Newtonian flow problem is presented; its finite-volume IBM approximation and the integration within the PIMPLE pressure-velocity coupling algorithm is introduced in Section \ref{sec:ibm}. Finally, some relevant examples of the application of the proposed method to industrial mixing processes are discussed in Section \ref{sec:results}.

\section{Motion of a generalized Newtonian fluid}  \label{sec:model}
The aim of this section is to introduce the differential equation model that describes the motion of a molten polymer when being processed. Molten polymers are often modelled by generalized Newtonian fluid that are special fluids whose viscosity is characterised by a nonlinear dependence on the shear rate and (possibly) temperature.
The differential model consists in the Navier-Stokes equations for non-isothermal incompressible fluids with shear and temperature dependent viscosity. 
In particular, we consider the mass, momentum and temperature conservation equations where the latter one is derived by the energy conservation principle under the hypothesis that the fluid is incompressible \cite{VERSTEEG}.

Let now $\vec{u}$ be the velocity, $p$ be the pressure and $T$ be the temperature fields. Then, the resulting system of equations that describes the motion of an incompressible generalized Newtonian fluid, reads:
\begin{equation}
\begin{aligned}
&\rho \pdd{\vec{u}}{}{t} +\rho  (\vec{u}\cdot \nabla) \vec{u}
- \nabla\cdot\left(\mu(\dot{\gamma},T)
(\nabla\vec{u} + \nabla^\intercal\vec{u})\right)
=-\nabla p
+ \vec{f},\\
&\nabla \cdot  \vec{u} = 0, \label{eq:goveq}\\
&\rho \pdd{c_p T}{}{t} + \rho\vec{u}\cdot\nabla (c_p T)
-\nabla\cdot(k \nabla T)=
\mu(\dot{\gamma},T)(\nabla\vec{u} + \nabla^\intercal\vec{u}):\nabla\vec{u} +
\rho r - \rho\vec{u}\cdot\vec{f},
\end{aligned}
\end{equation}
where $\mu$ and $\nu$ are the dynamic and kinematic viscosity, respectively, $c_p$ is the specific heat capacity, $k$ is the thermal conductivity. $\dot{\gamma}$ is the shear rate and $\vec{f}$ and $r$ are the momentum and energy source or sink terms, respectively.

Notice that, in the temperature equation, we consider also the contribution given by the viscous dissipation, which significantly contributes to the heating of the fluid, due to the high viscosity of polymers.

Equations \eqref{eq:goveq} are not closed until we define the constitutive relation for viscosity $\mu(\dot{\gamma},T)$. Viscosity is defined as the proportionality factor between stress and shear rate, hence the relation $\boldsymbol{\tau} = 2\mu\vec{D}$ holds for an incompressible fluid. Many different rheological laws can be considered to model the viscosity $\mu$. We consider the power law model, namely:
\begin{equation}
\mu(\dot{\gamma},T)=H(T)K\left(\dot{\gamma}\right)^{n-1},
\label{eq:power_law}
\end{equation}
where the parameters $K$ [Pa$\cdot$s] and $n$ represent the consistency factor, i.e. the zero-shear viscosity, and the power law exponent of the fluid, respectively. In this work, we will consider polymeric fluids with shear-thinning behaviour, corresponding to $n<1$.
The temperature shift factor $H(T)$ accounts for the possible dependence of viscosity on temperature, as it is typically the case for elastomers and thermoplastics. Here, we will consider the Arrhenius models, that reads:
\begin{equation}
\ln H(T) = \alpha \left(\cfrac{1}{T} - \cfrac{1}{T_r}\right),
\label{eq:arrhenius_law}
\end{equation}
where $\alpha$ is the activation temperature and $T_r$ [K] is the reference temperature.

\section{Non-conforming approximation of complex geometries}  \label{sec:ibm}
The aim of this section is to introduce the methods and algorithms that we used to solve the problem described in Section \ref{sec:model} and, in particular, we will set the focus on the review and formalization of the Immersed Boundary method (IBM), that has been used to approximate complex geometries on non-conforming grids.

The primary method that is employed to equations \eqref{eq:goveq} is the Finite Volume Method (FVM), specifically the one implemented in the open-source C++ library OpenFOAM (OpenFOAM.org version 10) \cite{WELLER1998, JASAK1996}, a CFD library widely adopted in industrial contexts. The FVM is a discretization technique which exploits the conservation of physical quantities. It uses an Eulerian approach, dividing the domain into control volumes and writing a local balance for each. By applying the Gauss-Green theorem, this becomes a balance of fluxes over the boundaries of the control volumes. These fluxes are then discretized using suitable numerical schemes. The method is also known for its adaptability with general polyhedral meshes and the ease with which it can be implemented in parallel to solve large-scale problems.


In addition to this, alongside to the FVM, the family of SIMPLE \cite{PATANKAR1980} methods is employed to solve the Navier-Stokes equations in a segregated manner. This technique, known as the projection method \cite{ELMAN2008}, allows for the separate resolution of momentum, temperature, and pressure equations, iteratively converging to the solution of the numerical problem.

On top of these methods, we introduce an IBM in order to deal with the complex geometries involved in polymer mixing processes and to avoid the generation of conforming grids, which can frequently prove to be prohibitively expensive in terms of computational resources. The IBM that we have implemented is a revised version of the one proposed in \cite{JASAKIBM2014} and implemented in \textit{foam-extend-4.0}, a fork of OpenFOAM. In particular, we will formalise the method and point out how we have improved the original method in terms of accuracy, robustness and parallel performances. We will also introduce a novel approach to integrate the non-conforming approach into segregated solution algorithms, like the SIMPLE, PISO and PIMPLE ones. This integration is crucial for enhancing the robustness of the Immersed Boundary Method (IBM) when applied to general polyhedral grids.

\subsection{The Immersed Boundary Method (IBM)}
We first define the FV mesh as a generic polyhedral tessellation of the physical domain.
Let $\domain\subset\Real^d$, where $d=3$, be a Lipschitz bounded domain. Let $\mesh$ be a polyhedral tessellation of $\domain$. Polyhedral cells of $\mesh$ are non-overlapping. $h=\max_{\cell_i\in\Tau_h}h_i$, where $h_i = \text{diam}(\cell_i)$ is the diameter of $\cell_i$. $\cellCentre_i$ and $\abs{\cell_i}$ denote the barycentre and the volume of cell $\cell_i$, respectively.

Let now $\hadomain$ be a hold-all domain such that $\domain\subset\hadomain$. Let $\Sigma$ be the surface triangulation of a closed manifold, representing the immersed boundary. Denote with $\triangle_i$ a triangle of $\Sigma$ and with $\abs{\triangle_i}$, $\triangleCentre_i$ the area and the barycentre of the $i$-th triangle, respectively. We denote with $\overline{\Sigma}$ the volume enclosed by the surface $\Sigma$.

Without loss of generality, suppose the flow to be external to the immersed surface. A generic Immersed Boundary method requires to ho through the following steps:
\begin{enumerate}
\item generate the triangulated surface $\Sigma$ representing the immersed object;
\item use the surface to divide the mesh in three regions: the \textit{solid set} $\Gamma_S$, containing the cells with barycentre inside the surface, representing the immersed object; the \textit{IB set} $\Gamma_\text{IB}$, containing the cells just outside the solid set, that is the  non-solid cells which have at least one solid cell as a neighbour; the \textit{fluid set} $\Gamma_F$, containing the remaining cells, representing the fluid domain;
\item generate the IB cell stencil $\S_i$, for each IB cell $\cell_i$, that is a set of mesh cells which are selected within certain criteria;
\item generate the IB approximator on each stencil (i.e. interpolation matrices);
\item at each time step and for each field, using the values of the actual solution evaluated on the stencil of an IB cell, evaluate the approximator function in the IB cell barycentre and impose the value on the solution.
\end{enumerate}

\begin{figure}[b]
\centering
\includegraphics[width=0.26\linewidth]{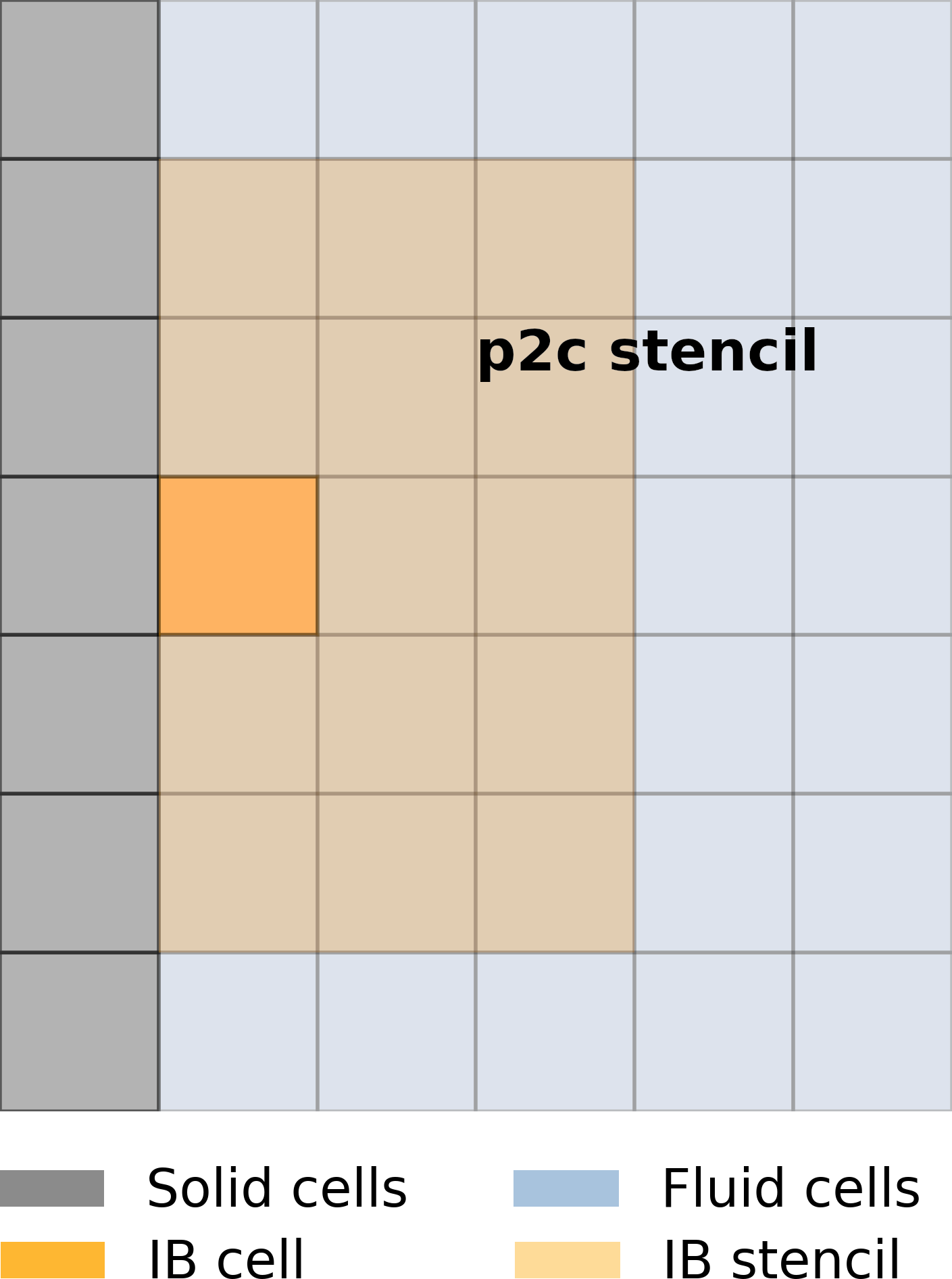}
\hspace{1cm}
\includegraphics[width=0.52\linewidth]{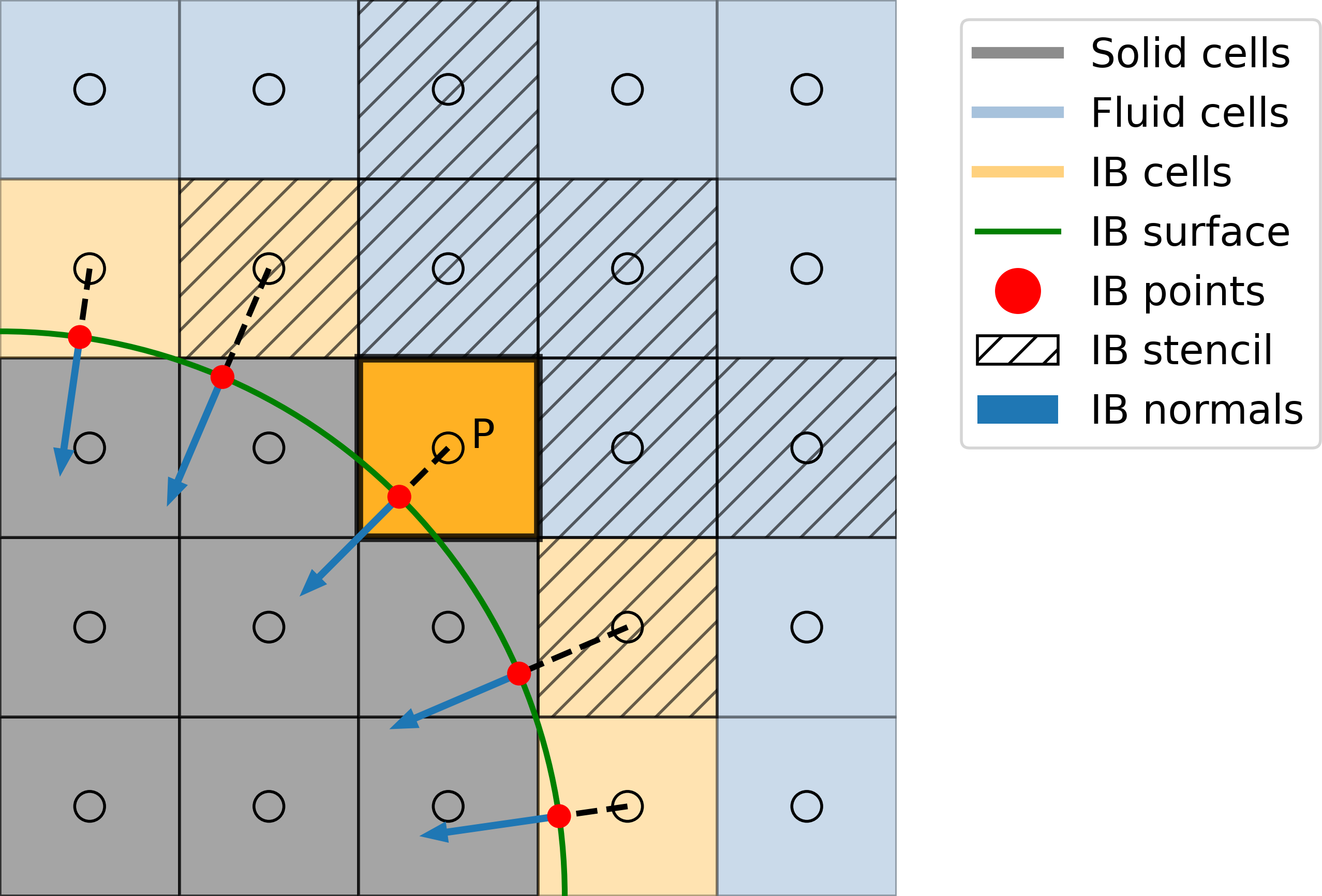}
\caption{Representation of point-to-cell stencil of level 2 (left) and schematics of IBM mesh elements (right).}
\label{fig:extended_stencil_filters:c2c_p2c_stencils}
\end{figure}

\begin{figure}[t]
\centering
\includegraphics[width=\linewidth]{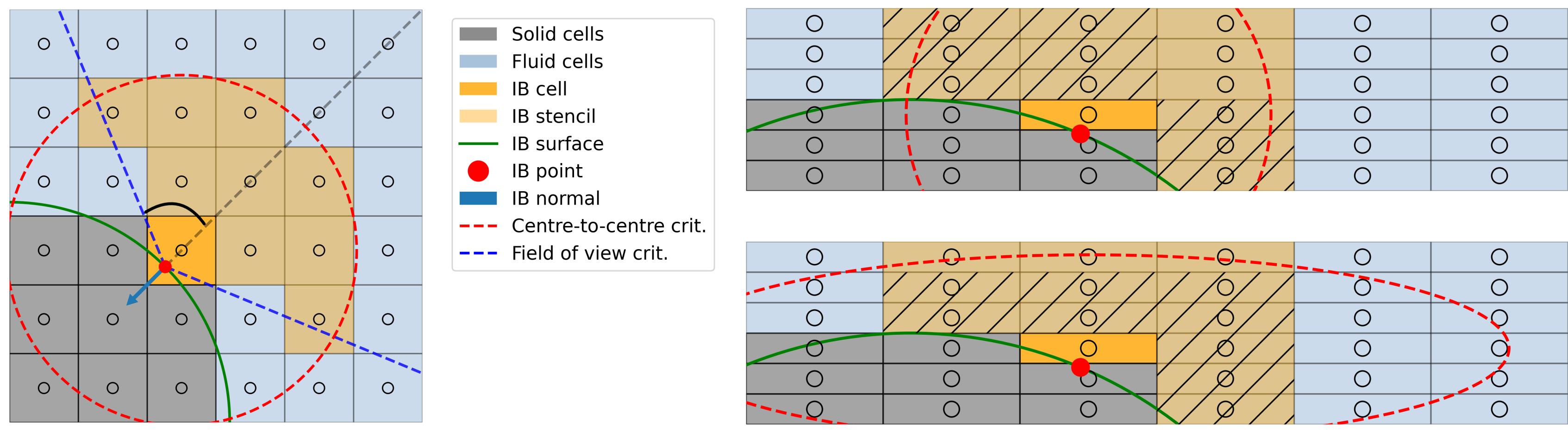}
\caption{Extended stencil of an IB cell (left), points cap criterion on anisotropic meshes (right): the closest 8 cells are hatched using standard Euclidean distance (top) and anisotropic distance (bottom).}\label{fig:extended_stencil_filters}
\end{figure}

For each IB cell $\cell_i\in\Gamma_{\text{IB}}$ we define the projections on $\Sigma$ of its cell center and we denote it as the \textit{IB points} $\pointvec_{\text{IB},i}$. Moreover, we denote with $\vec{n}_{\text{IB}}$ the normal unit vectors in the IB points directed inwards the solid set, as shown in Figure \ref{fig:extended_stencil_filters:c2c_p2c_stencils} (right).
The extended stencil $\S_i$ of the IB cell $\cell_i\in\Gamma_{\text{IB}}$ is constructed by collecting some neighbouring cells among IB and fluid sets. In particular, we choose cells whose barycentres simultaneously satisfy the following three criteria:
\begin{itemize}
\item \textit{spatial distance}: if centre-to-centre distance with respect to the IB cell is within a certain bound;
\item \textit{connectivity distance}: if the cell belongs to the stencil $\cells^c_i$ of level $c$ 
 (see Figure \ref{fig:extended_stencil_filters:c2c_p2c_stencils}, left) defined through the following recursive expression :
\begin{equation}
\cells^0_i = \{\cell_i\},\,\dots,\,
\cells^c_i = \{\cell_j \in \mesh: \cell_j \cap \cells^{c-1}_i \neq \emptyset\}\cup \bigcup\limits_{n=0}^{c-1}\cells^n_i.
\label{eq:extended_p2c_stencil_levelc}
\end{equation}
\item \textit{field of view}: if cell centre is within a certain angle with respect the the IB normal unit vector.
\end{itemize}
The combination of the first three criteria leads to the following definition of the extended stencil $\S_i$ of the IB cell $\cell_i\in\Gamma_{\text{IB}}$ (represented in Figure \ref{fig:extended_stencil_filters}, left):
\begin{equation*}
\S_i = \Bigg\lbrace K_j\not\in\Gamma_{S}: j\neq i,\cell_j \in \cells_i^{\bar{c}}, \norm{\cellCentre_i - \cellCentre_j}_2 \leq \bar{r},
    -\vec{n}_{i}\cdot\cfrac{(\cellCentre_i - \cellCentre_j)}{\norm{\cellCentre_i - \vec{p}_i}_2}\leq\cos(\bar{\theta}), \Bigg\rbrace,
\end{equation*}
where $\bar{c}$ is the maximum connectivity level, $\bar{r}$ is the limit distance and $\bar{\theta}$ is the field of view angle.

However, the application of these criteria frequently results in exceedingly large stencils. To address this, we introduced a new criterion that restricts the quantity of stencil points to a user-defined number, the \textit{points number cap}. In particular, we selects the $n$ cells in the stencil closest to the IB point according to an anisotropic distance metric weighted by the principal directions of the IB cell, namely:
\begin{equation}
\text{dist}(\point_{\text{IB}}, \point_i) = \norm{\Lambda^{-1}\mat{T}_{\text{PD}}(\pointvec_{\text{IB}}- \pointvec_i)},
\end{equation}
where $\mat{T}_{\text{PD}}$ is the principal direction tensor and $\Lambda$ is the diagonal matrix of its eigenvalues (see Figure \ref{fig:extended_stencil_filters}, right). 

Now we have to build a function approximator to impose the IB condition. First we define the \textit{IB value} $\vec{g}_{\text{IB}}=[g_i]$ as the vector of values of the immersed boundary condition evaluated at the IB points associated to IB cell $K_i\in\Gamma_{\text{IB}}$. We define the WLS approximant as a polynomial of degree $p$. Given  $\vec{x}=[x,y,z]^\intercal$, the polynomial basis $\polyspace_{\vec{x}} = [1,x,y,z,xy,\dots,y^{p-1}z,x^p,y^p,z^p]^\intercal$
and a vector of unknown coefficients $\boldsymbol{\beta}=[\beta_0,\dots ,\beta_{N_{coeffs}}]^\intercal$, we define the approximation function as
$f_{\text{IB}}(\vec{x}) = \polyspace_{\vec{x}}^\intercal \boldsymbol{\beta}.$

The WLS approximation is performed determining a set of observations of the solution, given by cell values of the extended stencil and the IB condition evaluated on the correspondent IB point, and computing the corresponding weights.

Let $u_j$ denote the FVM solution on the $j$-th cell of stencil $\S_i$. For each observation define also the weights $w_j \in [0,1],j=0,\dots,n$ associated to the observations points, where 0 is the IB point as a convention. We want to determine $\boldsymbol{\beta}^*$ such that:
\begin{equation*}
\boldsymbol{\beta}^*
=\argmin{\boldsymbol{\beta}}\norm{\text{W}^{\frac{1}{2}}\left(\vec{U}_{\text{IB}} - \text{A}\boldsymbol{\beta}\right)}^2
\label{eq:wls_problem}
\end{equation*}
where W is a square diagonal matrix with $\mat{W}_{jj}=w_j$, $\vec{U}_{\text{IB}}$ is the vector collecting all the of observations, and $A$ is the matrix that has on the $j$-th row the basis $\polyspace_{\vec{x}}$ evaluated in the $j$-th observation point, $\left[\text{A}\right]_{ij}=\left[\polyspace_{\vec{x}_i}\right]_j$. The weights have been defined as in \cite{JASAKIBM2014}. For the sake of simplicity we consider only Dirichlet condition on the immersed boundary (see \cite{NEGRINIPHD} for details on Neumann conditions).

Finally, let $S_{\text{IB}}$ be the interpolation operator that corrects the solution field on $\Gamma_{\text{IB}}$:
\begin{equation}
\vec{U}^{\text{corr}} = S_{\text{IB}}(\vec{g},\vec{U})=
s_{\text{IB}}g_{\text{IB},i} + \sum\limits_{\cell_j \in \S_i}s_j \vec{U}_j,
\label{eq:IB_interpolation_matrix}
\end{equation}
where $\vec{U}^{\text{corr}}$ is the corrected solution vector, $\vec{U}$ is the current solution vector, $\vec{s}_i=[s_{i,k}],k=1,\dots,\text{card}(\S_i)$ are the linear combination coefficients computed by $\vec{s}_i=\polyspace_{\cellCentre_i}^\intercal\boldsymbol{\beta}^*$. For what concerns the other regions, $S_{\text{IB}}$ extends the IB condition inside the immersed body on $\Gamma_{S}$ and is an identity on $\Gamma_F$. Also, it is worth noticing that, being $S_{\text{IB}}$ al linear combination, it can be represented by two matrices $\mat{S}_g$ and S, such that:
\begin{equation}
\vec{U}^{\text{corr}}=\mat{S}_g\vec{g} + \mat{S}\vec{U}.
\label{eq:ib_interp_matrix}
\end{equation}

\begin{remark}[Parallelization of IBM for large scale problems]  \label{rem:ibm:largescale}
In large-scale problems employing domain decomposition, IB cell stencils may cross processor boundaries. Without proper handling, this could lead to ill-posed WLS approximators due to insufficient interpolation points. To mitigate this, we use a communication map strategy, ensuring each processor exchanges only relevant data, thereby minimizing communication overhead. For each IB cell in current processor, we check which processors need to be queried to retrieve data related to its cell stencil. We gather all the ranks of these processors in a map that indicates where data should be sent and received. This optimizes the data communication and reduces the loss in computational performance with respect to the original implementation, as represented in Figure \ref{fig:scalability_old_ibm}.
\begin{figure}[t]
\centering
\includegraphics[width=0.9\linewidth]{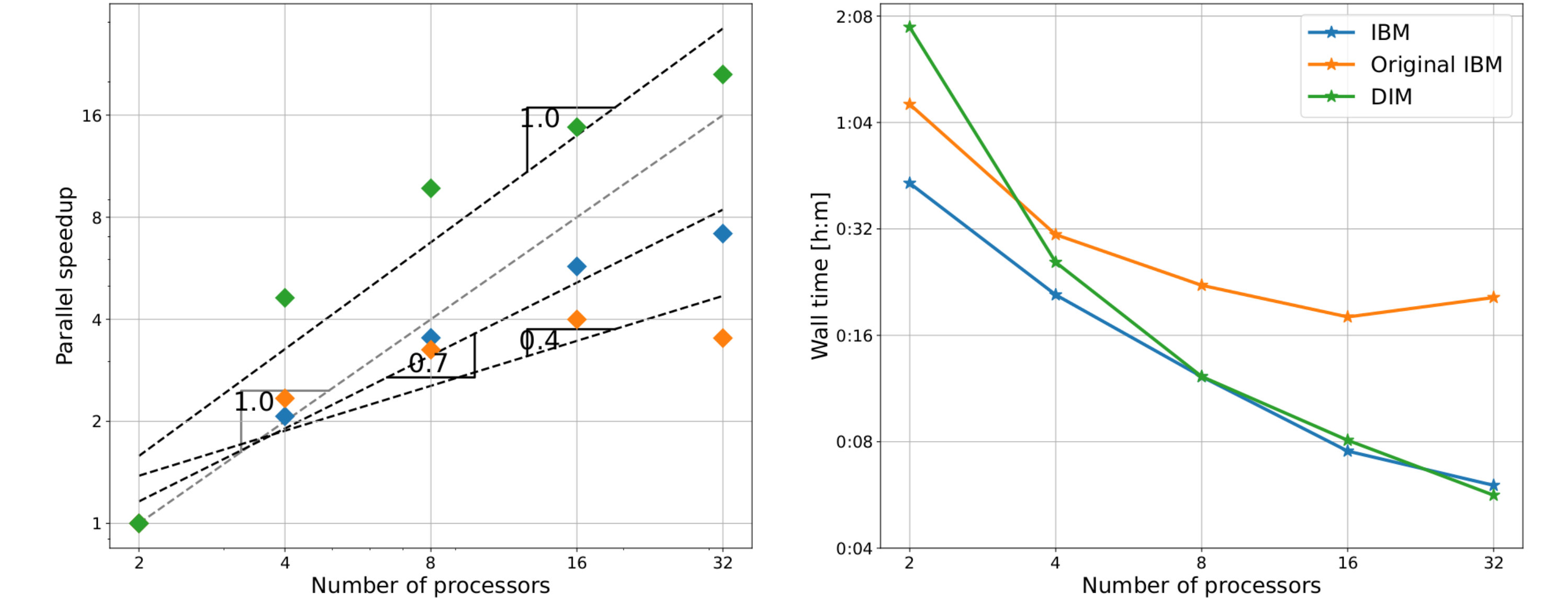}
\caption{Strong scalability test on a SSE test case, with around 4 million cells (cf. Section \ref{sec:ibm:numer}).}
\label{fig:scalability_old_ibm}
\end{figure}
\end{remark}

\begin{remark}[Enriched stencils]
While our discussion on stencil point selection has focused only on neighbouring cell centres, we have included also the possibility to consider boundary points from both conforming and immersed boundaries. This enhancement increases robustness in handling narrow boundary gaps, but also effectively manages multiple interacting immersed boundaries within the same simulation, such as screw-barrel and screw-screw gaps, and twin-screw extruder (refer to Section \ref{sec:results}).
\end{remark}

\begin{remark}[Diffuse Interface Method]  \label{rem:dim}
The approximator $f_{\text{IB}}$ is built such that any polynomial function can be used to correct the IB condition. In particular, setting $p=0$ we obtain the so called Diffuse Interface Method \cite{NEGRINI2021, SCHLOTTBOM2014, SCHLOTTBOM2016}.
\end{remark}

\subsection{Integration of IBM in the PIMPLE method}  \label{sec:ibm:sol}
As aforementioned, we use the SIMPLE family of projection methods to solve the Navier-Stokes equations. In general, projection methods require not only a decoupling of momentum and mass equations but also a splitting of velocity and pressure boundary conditions, hence a special treatment is needed for an immersed condition.

Here we consider the PIMPLE algorithm. We consider the time dependent isothermal incompressible Navier-Stokes equations for a Newtonian fluid complemented with suitable initial and boundary conditions. Given a time discretization scheme, the PIMPLE algorithm consists in three steps. For each time instant $t^{n+1}$, we first compute a velocity prediction, solving the momentum equation using pressure at time $t^n$. Then, we solve a Poisson equation for pressure, derived from the incompressibility constraint. Finally, we correct the velocity projecting it on divergence-free functions space, applying a pressure-dependent correction.

The procedure to integrate IBM and PIMPLE algorithm is based on imposing a consistent condition on pressure \cite{IKENO2007}. First, consider the algebraic counterpart (left) of Navier-Stokes system (right):
\begin{equation}
\begin{aligned}
\mat{A}\vec{U}^{n+1} + \mat{B}^\intercal\vec{P}^{n+1} &= \vec{b}, \qquad&&
\pdd{\vec{u}}{}{t} + (\vec{u}\cdot\nabla)\vec{u} -
\nabla\cdot\boldsymbol{\tau} + \nabla p = \vec{0}, \\
\mat{B}\vec{U}^{n+1} - \text{C}\vec{P}^{n+1} &= \vec{0}, \qquad&&
\nabla\cdot\vec{u} = 0
\end{aligned}
\label{eq:algebraic_ns}
\end{equation}
where $\mat{C}$ is the matrix that represents the Rhie-Chow stabilization term \cite{NEGRINI2023} and $\vec{b}$ is a source term used to impose the non-conforming boundary condition setting the IB and solid values in the region $\Gamma_{\text{IB}}\cup\Gamma_{S}$.

As we showed above, when dealing with non-conforming methods the velocity imposed in the IB region also depends on neighbouring values of the solution. Consider the IBM interpolation operator \eqref{eq:IB_interpolation_matrix}: $\vec{U}^*= S_{\text{IB}} (\vec{g},\vec{U}^*)$, where $*$ represents the algorithm step and $\vec{g}$ is the immersed boundary and solid region datum.

We now split the matrix A in its diagonal part D and its off-diagonal part -H, such that A=D-H. We also define matrix C as follows
\begin{equation}
\mat{C} = -\text{B}\text{D}^{-1}\text{B}^\intercal + \mat{R}(\text{D}^{-1})
\label{eq:rhie_chow_alge}
\end{equation}
where $\mat{R}(\text{D}^{-1})$ is the stiffness matrix with diffusivity coefficient $\text{D}^{-1}$ (see  \cite{NEGRINI2023}).

We divide the algorithm in five steps:
\begin{enumerate}
\item \textbf{Momentum predictor:} we make a prediction of the velocity field solving the first equation using pressure at step $n$, obtaining $\vec{U}^{n+\frac{1}{3}}$:
\begin{equation}
\text{A} \vec{U}^{n+\frac{1}{3}} + \text{B}^\intercal \vec{P}^{n}
= \vec{f} + \vec{b} \quad \implies \quad
\vec{U}^{n+\frac{1}{3}}=\text{A}^{-1}[- \text{B}^\intercal \vec{P}^{n}
+ \vec{f} + \vec{b}]
\label{eq:momentum_predictor}
\end{equation}

\item \textbf{Velocity correction:} we compute a new velocity field as
\begin{equation*}
\vec{U}^{n+\frac{2}{3}} = \text{D}^{-1}(\text{H} \vec{U}^{n+\frac{1}{3}} + \vec{f}).
\end{equation*}
that should satisfy the IB condition $\vec{U}^{n+\frac{2}{3}} =S_{\text{IB}}(\vec{g},\vec{U}^{n+\frac{2}{3}})$ in $\Gamma_{\text{IB}}\cup\Gamma_{S}$. Then, velocity field at step $n+1$ would read:
\begin{equation}
\vec{U}^{n+1}
= \vec{U}^{n+\frac{2}{3}} - \text{D}^{-1}\text{B}^\intercal \vec{P}^{n+1}.
\label{eq:momentum_correction}
\end{equation}
However, $\vec{P}^{n+1}$ is unknown at this step.

\item \textbf{Pressure equation:} by multiplying equation \eqref{eq:momentum_correction} by divergence matrix $\text{B}$, by definition we obtain $\text{B}\vec{U}^{n+\frac{2}{3}} - \text{C}\vec{P}^{n+1}= 0$, hence the pressure equation reads:
\begin{equation}
\text{B}\vec{U}^{n+\frac{2}{3}}
= \text{B}\text{D}^{-1}\text{B}^\intercal \vec{P}^{n+1}
+ \text{C} \vec{P}^{n+1}
= \text{R}(\text{D}^{-1})\vec{P}^{n+1}.
\label{eq:pressure_equation}
\end{equation}

\textbf{IB interpolation:} in general, when applying the IB interpolator $S_{\text{IB}}$ we are performing a matrix-vector multiplication (see equation \eqref{eq:ib_interp_matrix}). If we are in the DIM case (cf. Remark \ref{rem:dim}), the PIMPLE algorithm is not modified because $S_{\text{IB}}$ is an identity. In general, when the IB values are corrected, the pressure equation requires a modification.

Requiring that $\vec{U}^{n+1}$, satisfying equation \eqref{eq:momentum_correction}, is consistent with the non-conforming condition, we can derive a correction for the pressure field $\vec{P}^{n+\frac{2}{3}}$. The constraint reads: let $\mat{S}_g, \mat{S}$ be the matrices that represent the IB interpolator $S_{\text{IB}}$,
\begin{equation*}
\vec{U}^{n+1}
= S_{\text{IB}}(\vec{g},\vec{U}^{n+1})
= \mat{S}_g\vec{g} + \mat{S}\vec{U}^{n+\frac{2}{3}} - \mat{S}\text{D}^{-1}\text{B}^\intercal \vec{P}^{n+1}
\end{equation*}

By applying the latter correction, by linearity of the operators and by choosing opportunely the Rhie-Chow stabilization matrix, the pressure equation will read:
\begin{equation}
\text{B}\vec{U}^{n+1}
= \text{B}\mat{S}\text{D}^{-1}\text{B}^\intercal \vec{P}^{n+1}
+ \text{C} \vec{P}^{n+1} 
= \mat{R}(\mat{S}\text{D}^{-1})\vec{P}^{n+1}.
\label{eq:modified_pressure_equation}
\end{equation}

\item \textbf{Final correction:} compute $\vec{U}^{n+1}$ using equation \eqref{eq:momentum_correction}.
\end{enumerate}

This approach embodies a consistency property in our IBM when using a projection algorithm. Moreover, in our experience, this approach also boosts the robustness of the IBM in many scenarios with respect to the original implementation, where a homogeneous Neumann condition was imposed on pressure. In particular, we mention the case of anisotropic meshes, such as wall boundary layers, which are needed in order to simulate realistic mixing devices geometries avoiding an excessive number of degrees of freedom and to properly capture the flow field inside the gap between internal rotors and the external barrel.

\section{Numerical results on industrial continuous mixing devices} \label{sec:results}
In this section, we present a set of numerical results obtained with the finite-volume Immersed Boundary scheme introduced above applied to a family of continuous mixing devices of industrial interest including a Single-Screw Extruder (SSE), Twin-Screw Extruder (TSE) and a Planetary Rolling Extruder (PRE).
The different devices are characterized by different levels of geometrical complexity, as well as specific configuration features that require ad-hoc simulation settings, as it will be detailed in the following sections.
The objective of this section is to show the wide spectrum of applicability of the proposed methodology in the context of continuous mixing, considering real industrial geometries.

\subsection{Single-Screw Extruder}\label{sec:ibm:numer}
\begin{figure}[t]
\centering
\includegraphics[width=0.7\linewidth]{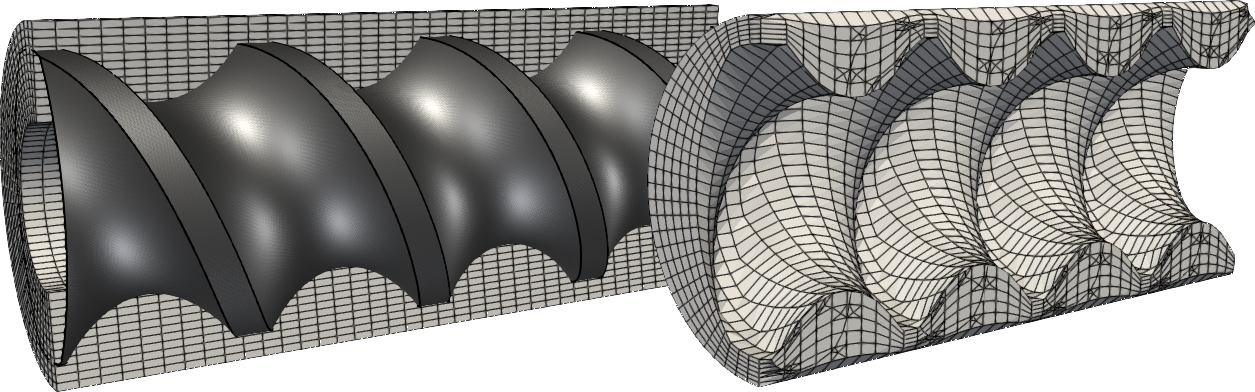}
\caption{Representation of the meshes employed to perform simulations. Left: uniform non-conforming mesh. Right: conforming mesh.}
\label{fig:ibm_conv:ssegeometry}
\end{figure}

This first case has been carried out to validate the accuracy of the IBM on a case of interest: the non-isothermal flow of a non-Newtonian fluid in a single-screw extruder, where the problem setup is the one that was proposed in \cite{KIM2006}. The geometry is represented in Figure \ref{fig:ibm_conv:ssegeometry} and is built by twisting, i.e. extruding and rotating, the section along the axial direction. The barrel is represented by a cylinder. The conforming mesh has been generated by deforming an hollow cylinder, moving the inner cylinder points to conform the screw surface and the internal points accordingly.
For the boundary conditions, a rotational velocity of 90 RPM is imposed on the screw and a no-slip condition on the barrel. Pressure conditions are set accordingly. We then imposed a temperature of 473K on the barrel while the inflow temperature is set to a uniform distribution at 463K. For the screw we considered an homogeneous Neumann condition.

\begin{figure}[h]
    \includegraphics[width=\linewidth]{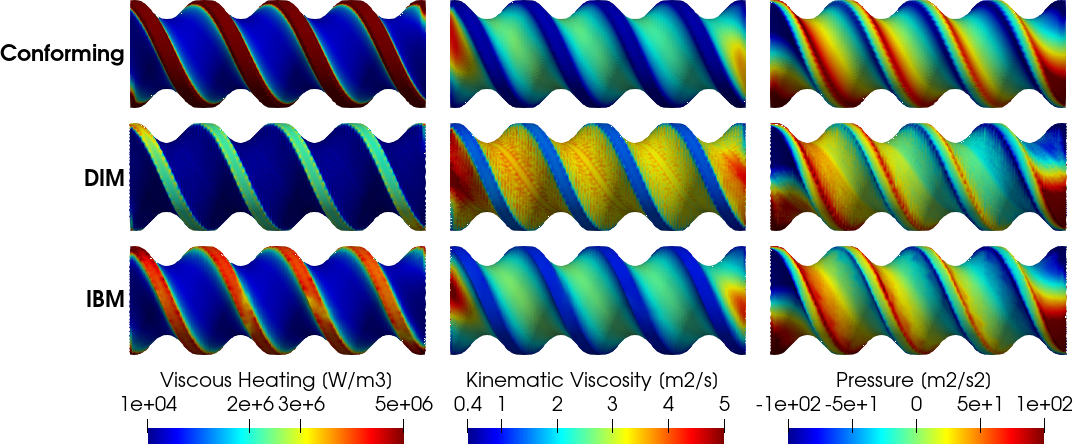}
\caption{Evaluation of various quantities on the screw surface for different discretization methods.}
\label{fig:ibm_conv:sse3d_surf_neu}
\end{figure}

We consider three methods to approximate the geometry: the geometry-conforming method, the DIM and the IBM, where we take the conforming solution as reference. In particular, we compare some quantities evaluated on the screw surface using an extrapolation based on the IB approximator (see Section \ref{sec:ibm:numer}). The results are reported in Figure \ref{fig:ibm_conv:sse3d_surf_neu}. For each quantity, the IBM profile is smoother than DIM profiles and it is much closer to the conforming one.

\subsection{Twin-Screw Extruder}
\begin{figure}[h]
    \centering
    \includegraphics[width=0.4\linewidth]{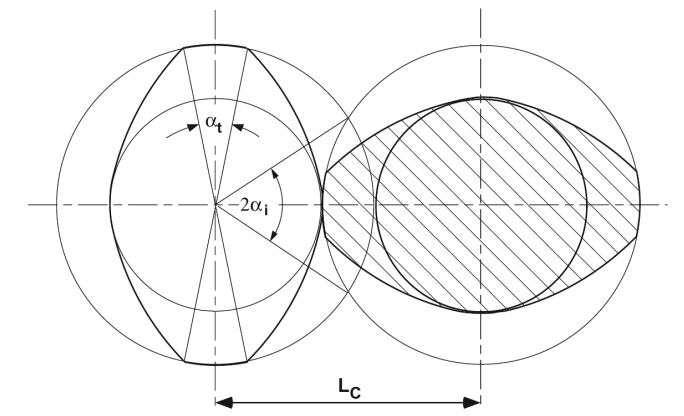}
    \hspace{10mm}
    \includegraphics[width=0.4\linewidth]{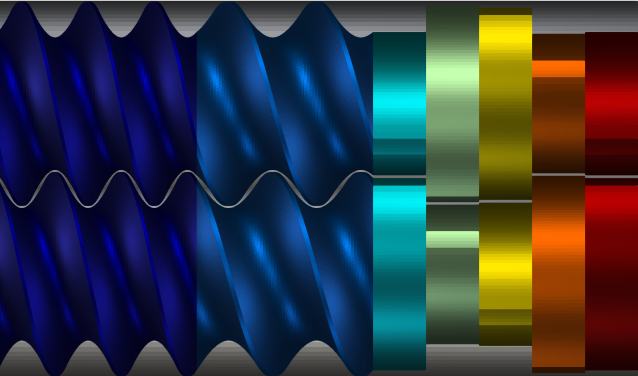}
    \caption{Sketch of the TSE section \cite{RAUWENDAAL2014} (left) and modular axial assembly (right).}
    \label{fig:tsegeom}
\end{figure}

The second mixing device that we consider is a twin-screw extruder (TSE), a well established technology extensively used for mixing, compounding, or reacting polymeric materials.  Larger heat transfer area and better mixing ability allow good control of stock temperatures, residence times and positive conveying, that are key elements in the extrusion of thermally sensitive materials. On the other hand, its geometry complexity makes the TSE more difficult to be simulated when compared to the SSE. These devices are typically designed with removable elements, in order to allow arbitrary sequences
of elements along the shaft. This modular design is highly flexible and allows process optimization. In this respect, numerical modelling may play a fundamental role in designing custom elements patterns for specific applications. We consider in this analysis a self-wiping intermeshing co-rotating TSE \cite{RAUWENDAAL2014}.
The section of a TSE is sketched in Figure \ref{fig:tsegeom}, left.
The device is composed by three parts: the barrel, represented by two adjacent cylinders of total width $2L_C$, of the order of dozens of centimetres, and length 30$L_C$ and the two screws positioned in the two channels. The
gap between screws and barrel is of the order of tens of millimetres. As mentioned before, the screws are composed by different elements. In our case we consider two transport modules with different pitches, five kneading modules and another transport module (see Figure \ref{fig:tsegeom}, right). The transport modules are devoted to conveying the
fluid towards the die with different velocities depending on the flight angle, however, their
mixing action is limited. The mixing action is prevalently performed by kneading modules,
that are basically screw elements with a 90 degrees flight angle. They can be also of different lengths: long modules for distributive mixing, while short ones for dispersive
mixing. 

\begin{figure}[b]
    \centering
    \includegraphics[width=0.325\linewidth]{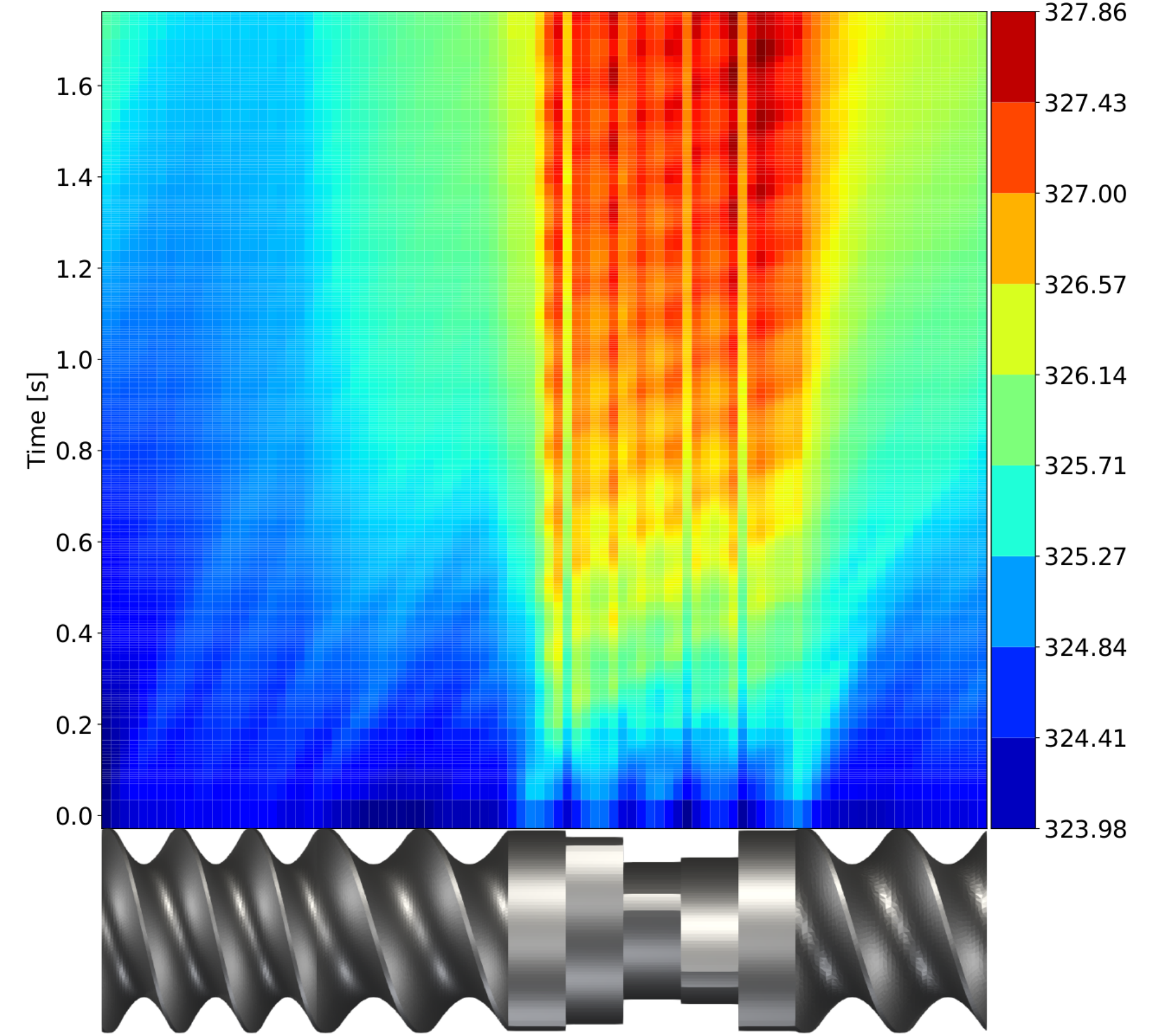}
    \includegraphics[width=0.325\linewidth]{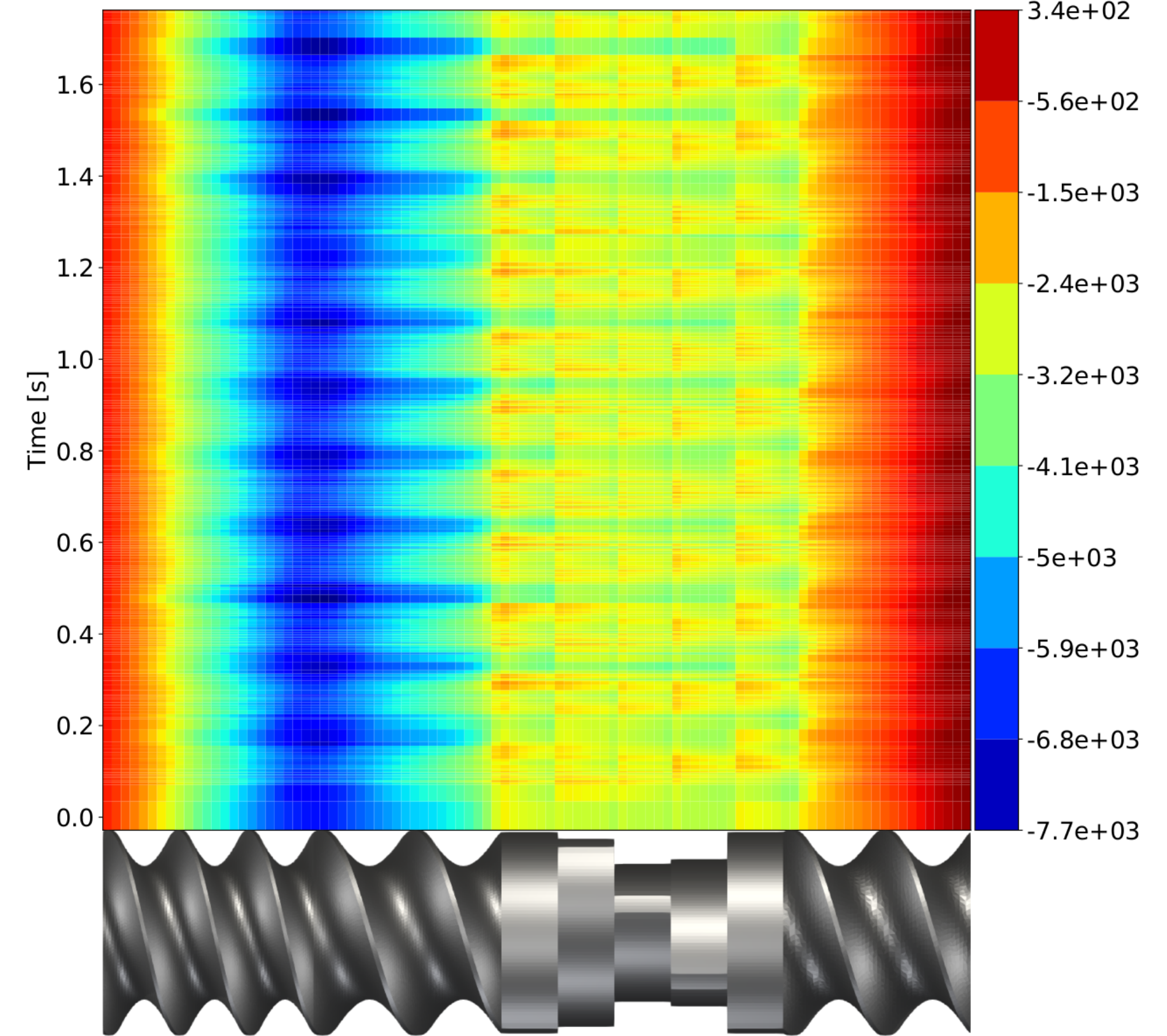}
    \includegraphics[width=0.325\linewidth]{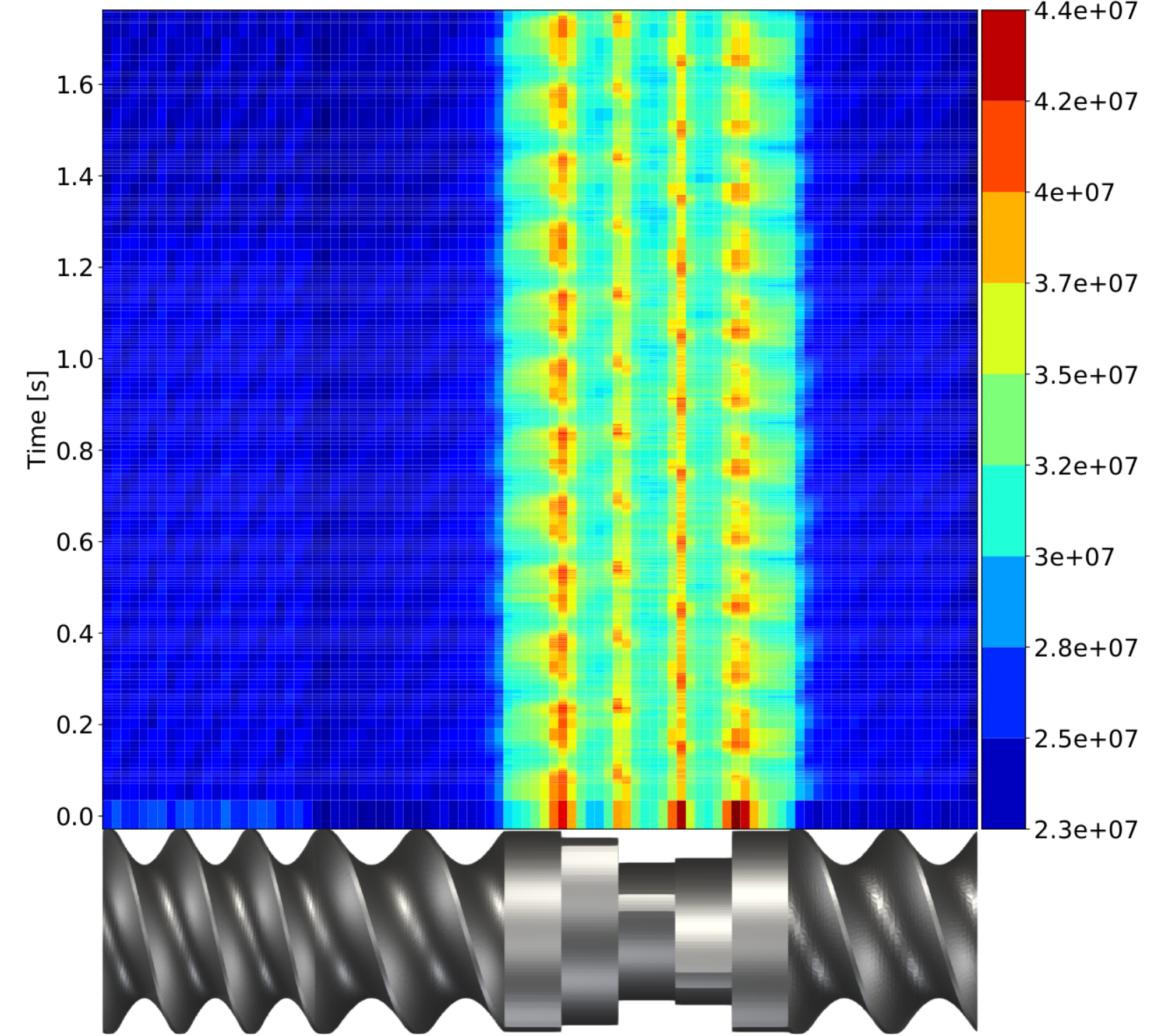}
    \caption{Time evolution of axial distribution of section averaged quantities: temperature average (left), pressure average (middle), and viscous heating average (right.}
    \label{fig:tsetransient}
\end{figure}

We present a time dependent simulation of the revolution of a TSE filled with a polymer with a high-viscosity power law rheology.
The working conditions are the following. The flow is driven by the screw rotation, set at 100RPM, and atmospheric pressure boundary conditions on both inflow and outflow. On barrel and screws, we set no slip condition for velocity. The barrel and screws temperatures are set to 323K as the inflow temperature. Three complete screw revolutions have been simulated and the transient behaviour of the flow and the temperature field has been analyse. As shown in Figure \ref{fig:tsetransient}, where the time evolution of the axial distribution of section averaged quantities is displayed, flow rate, pressure and viscous heating reach their asymptotic solution quite fast, while the temperature keeps increasing in time and would require a longer simulation time to reach convergence.

\begin{figure}[h]
    \centering
     \includegraphics[width=0.325\linewidth]{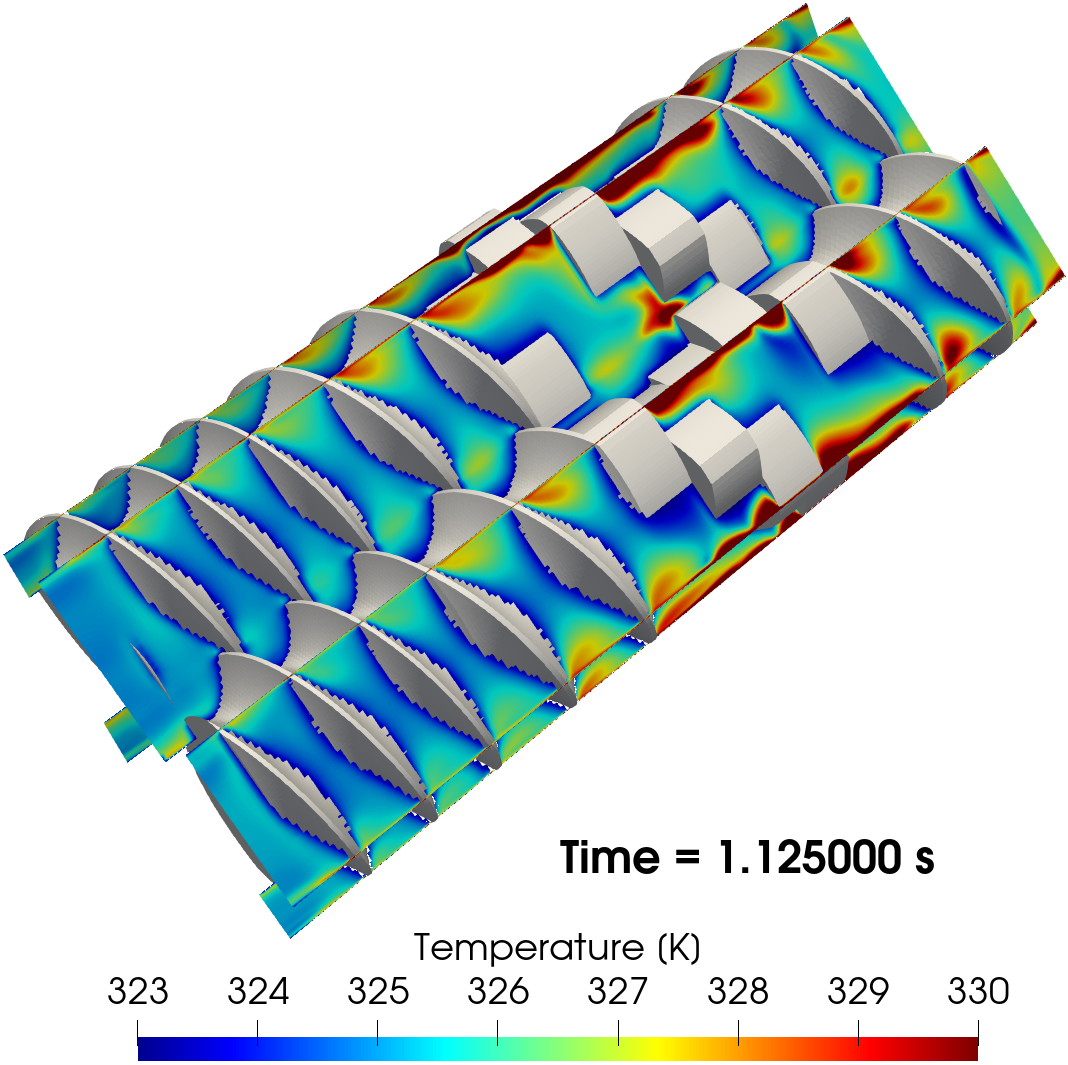}
     \includegraphics[width=0.325\linewidth]{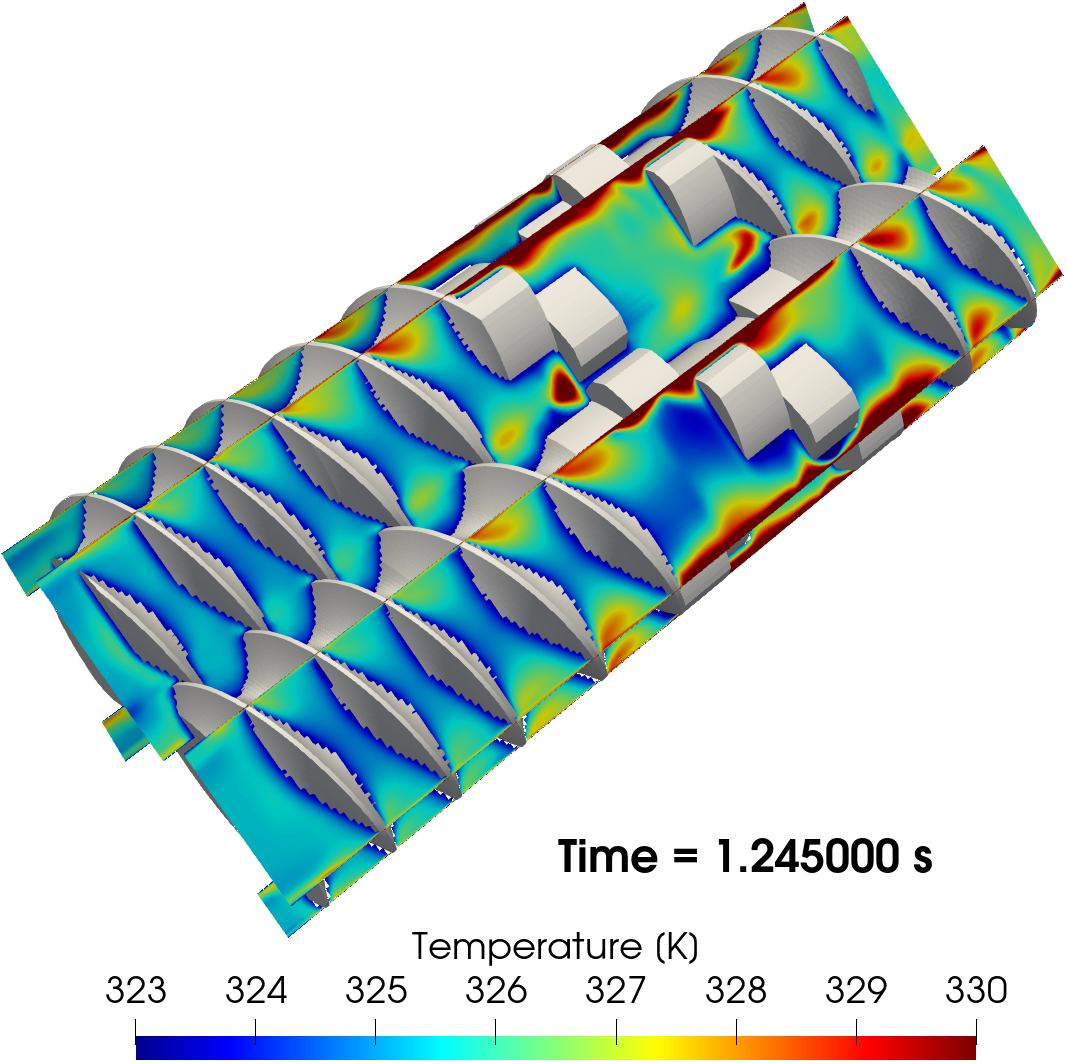}
     \includegraphics[width=0.325\linewidth]{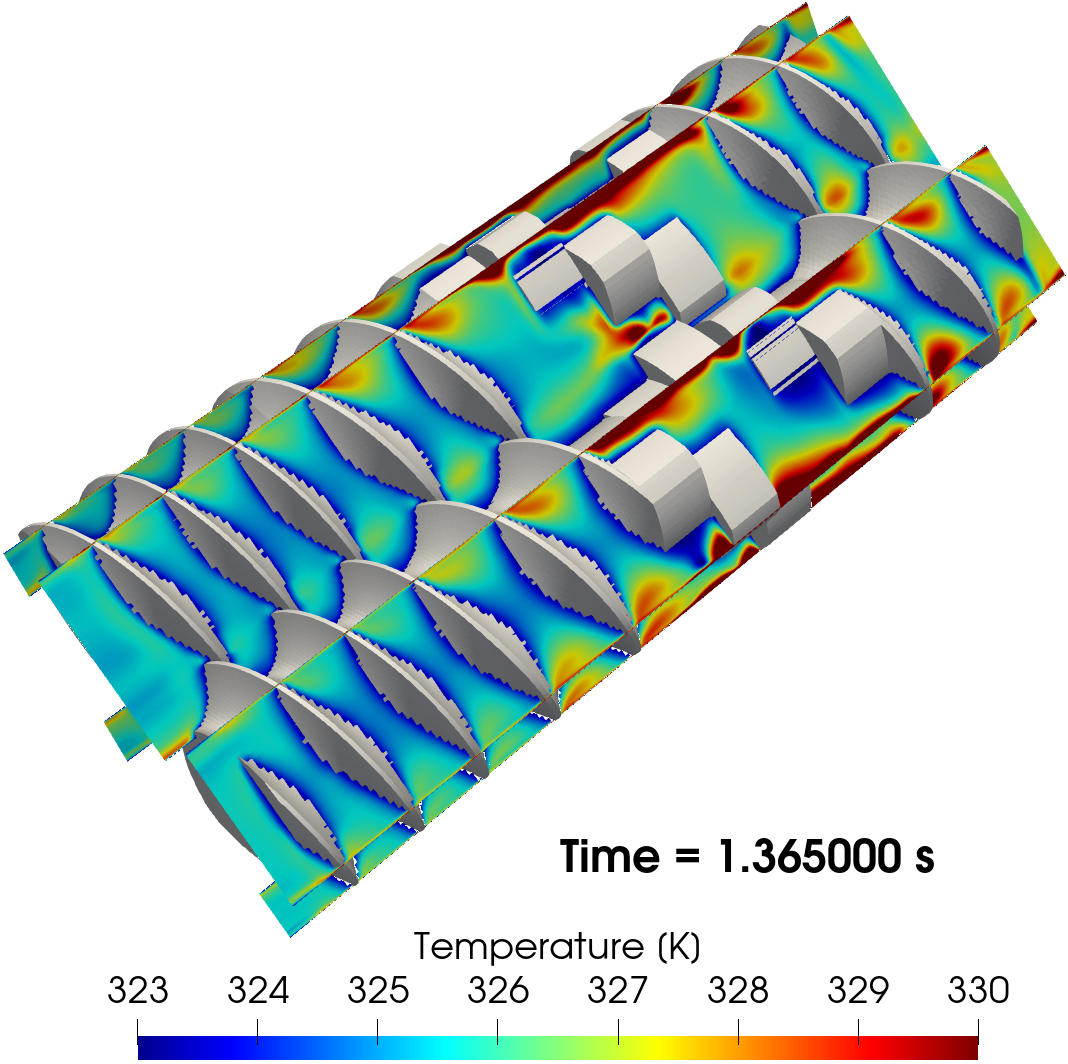}
    \caption{Temperature distribution in the twin-screw extruder at different time instants.}
    \label{fig:tsetemp}
\end{figure}

Temperature distribution inside the TSE at different time instants is displayed in Figure \ref{fig:tsetemp}. The small gap sizes, in combination with high rotational velocities of the screws, heat the fluid through high values of viscous heating. Indeed the highest temperature levels are reached in the gaps between the two screws, in particular between kneading modules.

\subsection{Planetary Roller Extruder}
\begin{figure}[h]
\centering
\includegraphics[width=0.3\linewidth]{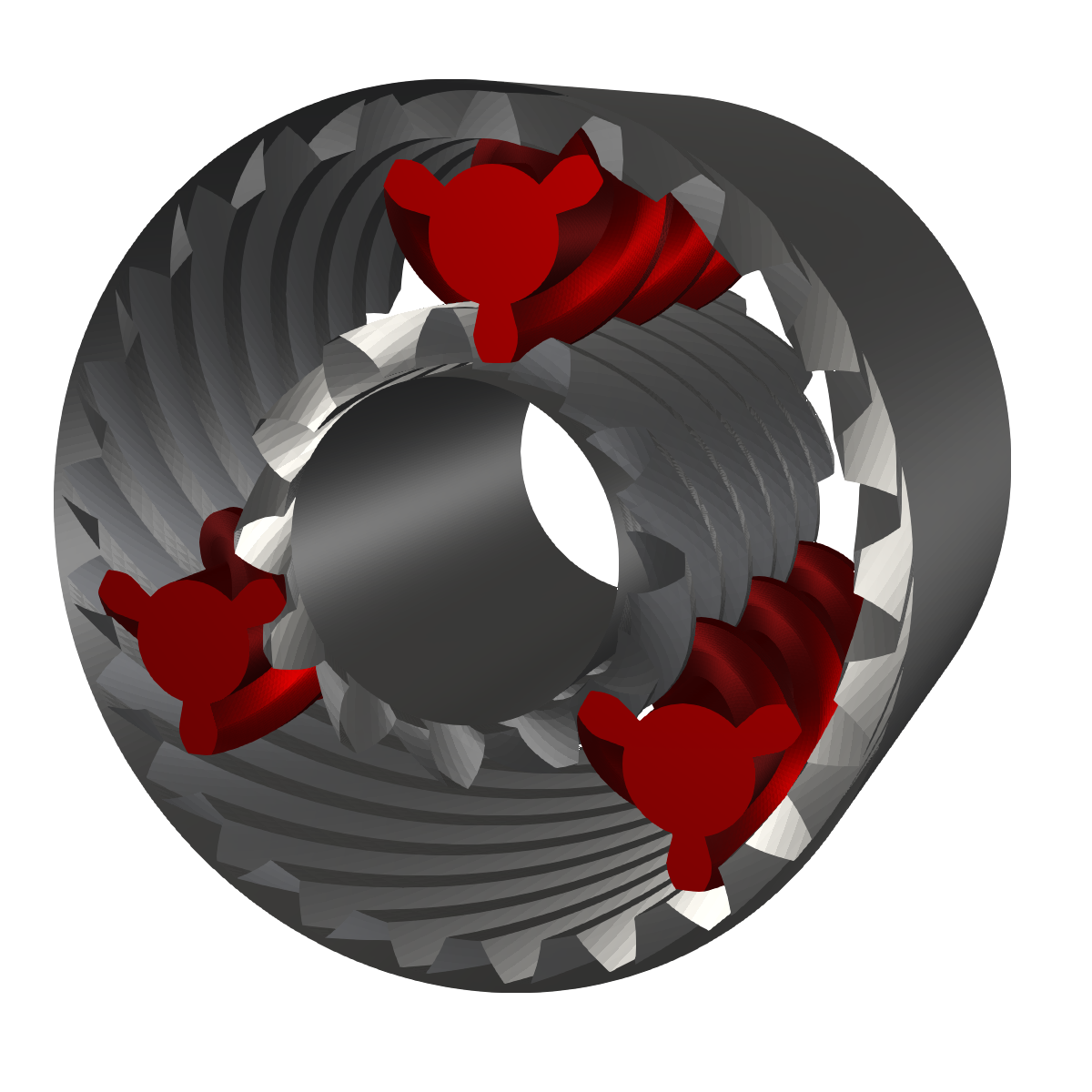}
\hspace{10mm}
\includegraphics[width=0.3\linewidth]{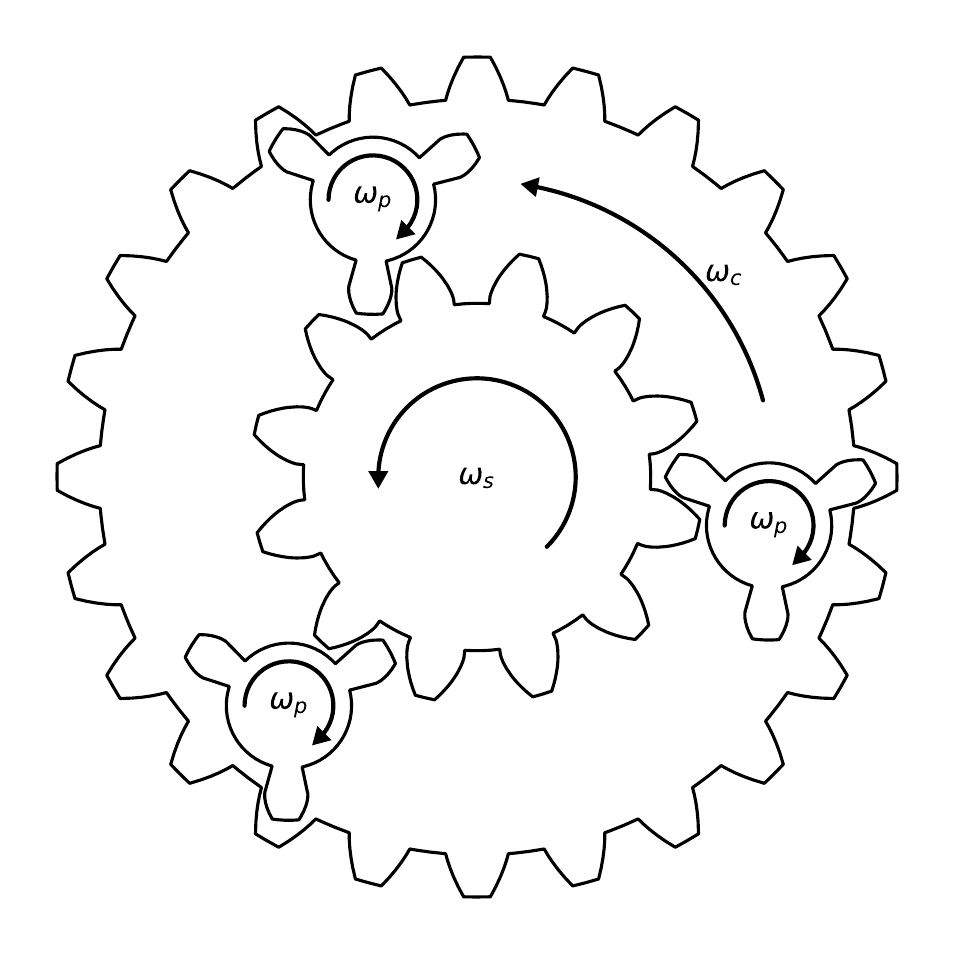}
\caption{PRE configuration (left) and 2D kinematic sketch (right).}
\label{fig:pregeom}
\end{figure}

In this last section, we present the results obtained using the Immersed Boundary method to simulate the flow of a polymeric material in a PRE device with 3 spindles (see Figure \ref{fig:pregeom}, left).
Given the radii of sun ($R_s$), ring ($R_r$) and planets ($R_p$) and considering a fixed ring, the kinematics is completely defined by the sun shaft rotation velocity ($\omega_r$), see Figure \ref{fig:pregeom}, right. In particular, the barycentres of the planets rotate around the
origin with velocity $\omega_c=\omega_s R_s/(2(R_s+R_p))$ while the planets rotates around their own barycentre with velocity $\omega_p=-\omega_c(R_s+R_p)/R_p$.

To simplify the implementation of this kinematics, the problem is solved on a non-inertial reference frame moving with the planet barycentres. In this reference frame each component rotates around its own axis with angular velocities:
\begin{equation}
\tilde{\omega}_s = \omega_s - \omega_c, \quad \quad
\tilde{\omega}_r = -\omega_c, \quad \quad
\tilde{\omega}_p = \omega_p - \omega_c.
\end{equation}

In order to minimize the accuracy limitations associated to non-conforming geometries, differently to what done in \cite{WINCK2021} where the mesh superposition technique was adopted for all boundaries, here only the spindles are treated as immersed boundaries, while the grid is conforming to sun and ring geometries.

\begin{figure}[t]
\centering
\includegraphics[width=0.35\linewidth]{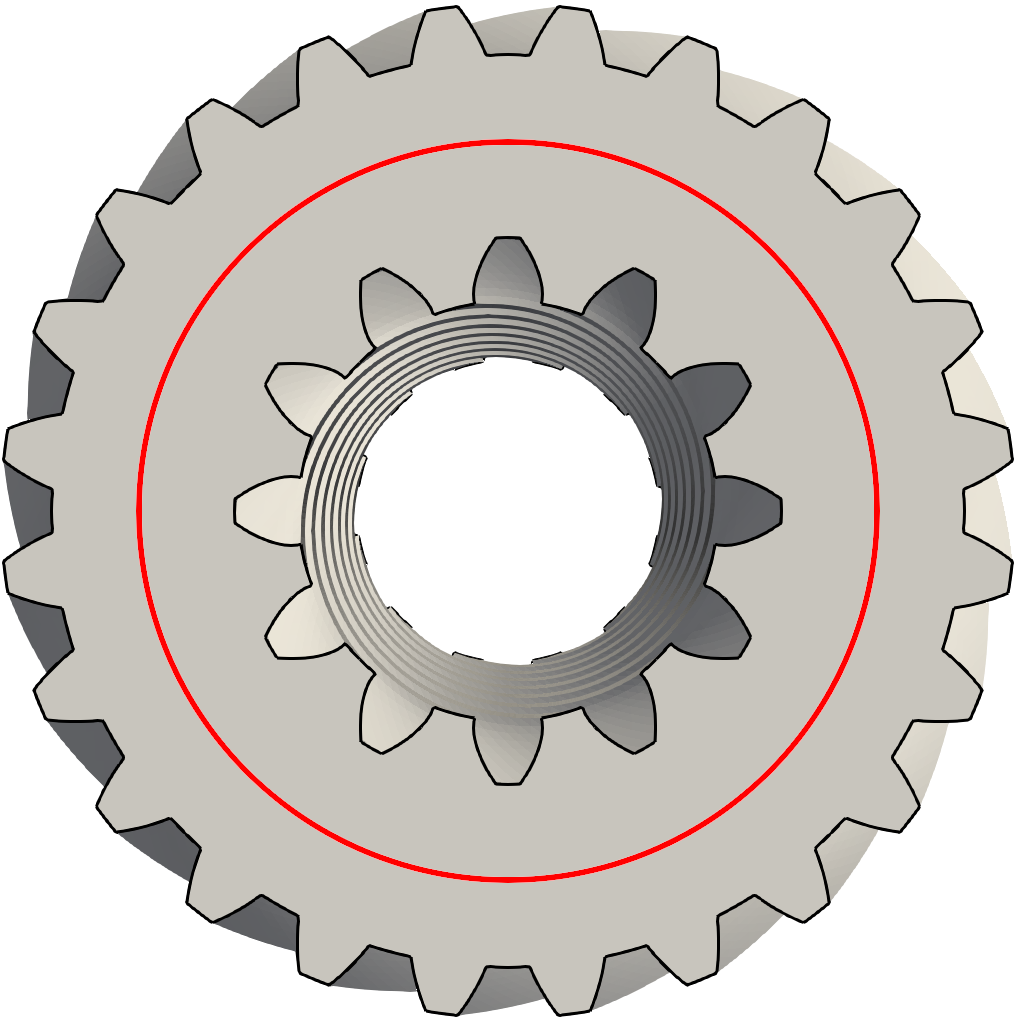}
\hspace{10mm}
\includegraphics[width=0.35
\linewidth]{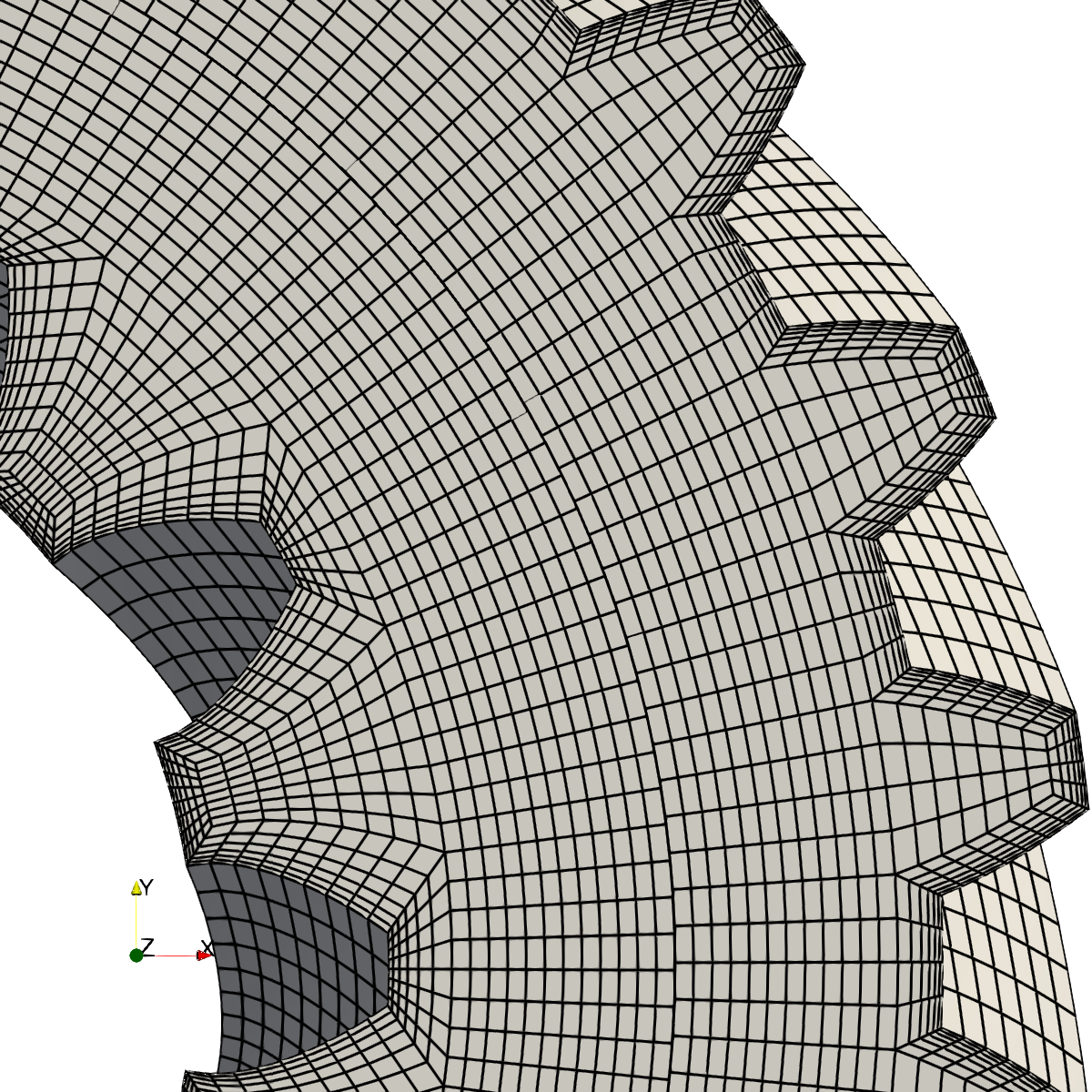}
\caption{Sub-domains separated by the sliding mesh interface (left) and non-conformal mesh motion (right).}
\label{fig:sliding}
\end{figure}

Achieving this result is not trivial since sun and ring have different angular velocity and a conforming grid treatment of sun and ring can only be obtained including a sliding mesh interface (in red in Figure \ref{fig:sliding}, left) that is dealt in OpenFOAM through the non-conformal coupling (NCC) conditions. During the motion the inner and outer portion of the mesh rotate in opposite direction and the continuity of fluxes is imposed on the non-conformal sliding interface (see  \ref{fig:sliding}, right). 

The flow rate and the mean temperature are imposed at the inflow while zero-mean pressure and homogeneous Neumann conditions on temperature are imposed at the outflow. No-slip condition and fixed temperatures are imposed on sun and ring. On the non-conforming planet surface IB no-slip conditions are imposed together with an interface condition prescribing the continuity of heat flux between the fluid and the solid. For additional details on the case setup and a comprehensive set of numerical results, we refer to \cite{NEGRINIPHD}.

\begin{figure}[h]
    \includegraphics[height=3.55cm]{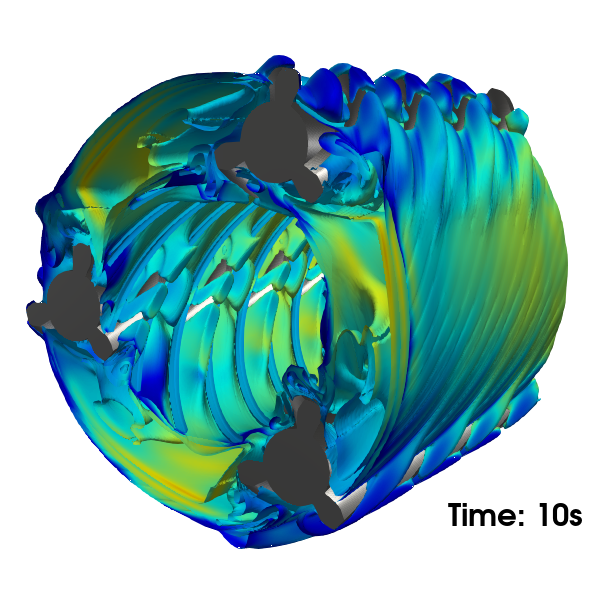}
    \includegraphics[height=3.55cm]{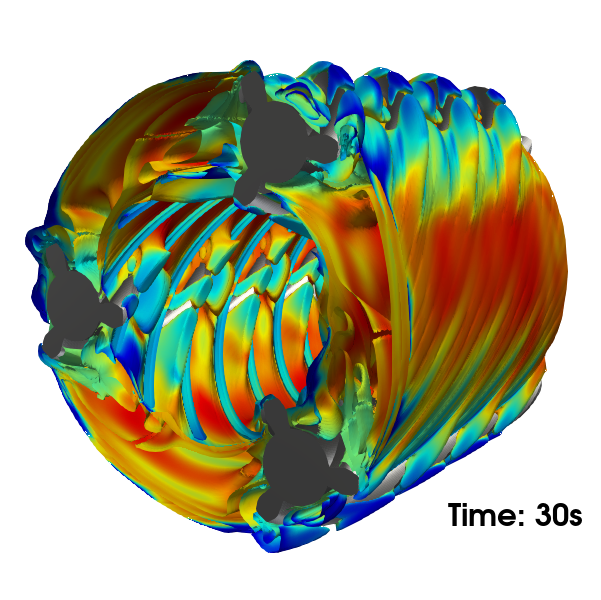}
    \includegraphics[height=3.55cm]{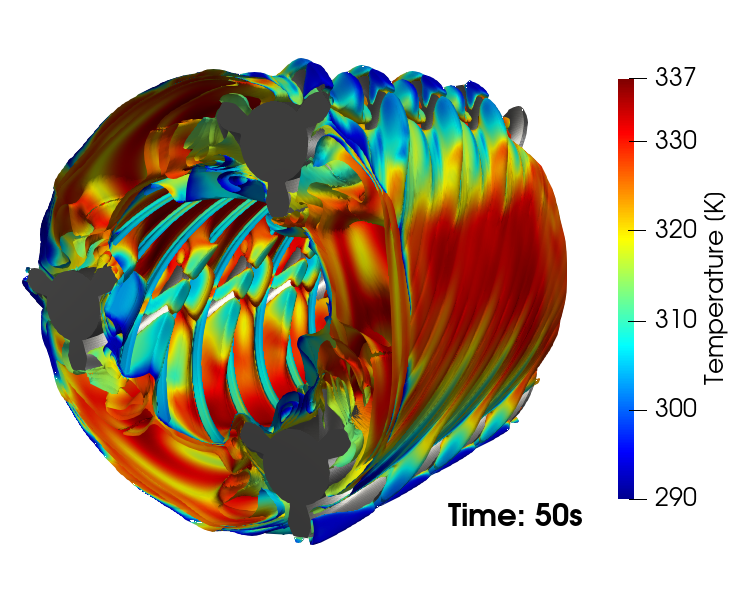}
\caption{Temperature distribution evaluated for different time instants on relative velocity magnitude iso-surfaces.}
\label{fig:pre}
\end{figure}

As an example illustrating the capability of the proposed method to tackle this complex configuration, in Figure \ref{fig:pre} the temperature distribution inside the PRE at different time instants is shown.

\section{Conclusions}
In this work, a new implementation of the finite-volume Immersed Boundary method in the OpenFOAM library has been introduced. The main goal of this numerical tool is the simulation of a range of polymeric mixing devices characterized by complex geometries and kinematics that cannot be faced by body conforming approaches. The accuracy and efficiency of the method has been investigated through a large set of simulations attesting that the presented method represents a valuable tool for mixing task analyses.

\section{Acknowledgements}
This research is part of the activities of Dipartimento di Eccellenza 2023-2027. The present work has been carried out in collaboration with the industrial partner Pirelli Tyre S.P.A., in the framework of the Pirelli-Polimi Joint Lab with the financial support of Fondazione Politecnico. N.P. and M.V. have been partially funded by PRIN2020 n.20204LN5N5 ``Advanced polyhedral discretisations of heterogeneous PDEs for multiphysics problems''. N.P. has been partially supported by ICSC–Centro Nazionale di Ricerca in High Performance Computing, Big Data, and Quantum Computing funded by the European Union–NextGenerationEU plan. M.V. is member of INdAM-GNCS. 

\bibliographystyle{ieeetr}
\bibliography{IBM_M2P_NPV}

\begin{thebibliography}{10}

\bibitem{JI2020}
D.~Ji, Y.~Xiao, Q.~Huang, and H.~Shi, ``Safety design and numerical simulation
  of twin screw extruder for energetic materials,'' {\em J.~Phys.~Conf.~Ser.},
  vol.~1507, p.~022027, mar 2020.

\bibitem{WINCK2021}
J.~Winck and S.~Frerich, ``Numerical simulation of fluid flow and mixing
  dynamics inside planetary roller extruders,'' {\em Int. Polym. Process.},
  vol.~36, no.~5, pp.~508--518, 2021.

\bibitem{MST}
X.~Zhang, L.-F. Feng, W.-X. Chen, and G.-H. Hu, ``Numerical simulation and
  experimental validation of mixing performance of kneading discs in a twin
  screw extruder,'' {\em Polym.~Eng.~Sci.}, vol.~49, pp.~1772 -- 1783, 09 2009.

\bibitem{SIRJALA2000}
S.~Syrjala, ``Numerical simulation of nonisothermal flow of polymer melt in a
  single-screw extruder: A validation study,'' {\em Numer.~Heat Transf.; A:
  Appl.}, vol.~37, no.~8, pp.~897--915, 2000.

\bibitem{CHANG1995}
R.-Y. Chang and K.-J. Lin, ``The hybrid {FEM/FDM} computer model for analysis
  of the metering section of a single-screw extruder,'' {\em Polym.~Eng.~Sci.},
  vol.~35, pp.~1748+, Nov 1995.

\bibitem{HELMIG2019}
J.~Helmig, M.~Behr, and S.~Elgeti, ``Boundary-conforming finite element methods
  for twin-screw extruders: Unsteady - temperature-dependent - non-newtonian
  simulations,'' {\em Computers {\&} Fluids}, vol.~190, pp.~322--336, 2019.

\bibitem{HINZ2020}
J.~Hinz, J.~Helmig, M.~M\"{o}ller, and S.~Elgeti, ``Boundary-conforming finite
  element methods for twin-screw extruders using spline-based parameterization
  techniques,'' {\em Comput. Methods Appl. Mech. Eng.}, vol.~361, p.~112740,
  2020.

\bibitem{PESKIN2002}
C.~S. Peskin, ``The immersed boundary method,'' {\em Acta Numerica}, vol.~11,
  pp.~479--517, 2002.

\bibitem{MITTAL2005}
R.~Mittal and G.~Iaccarino, ``Immersed boundary methods,'' {\em Annual Review
  of Fluid Mechanics}, vol.~37, no.~1, pp.~239--261, 2005.

\bibitem{BOFFI2003}
D.~Boffi and L.~Gastaldi, ``A finite element approach for the immersed boundary
  method,'' vol.~81, pp.~491--501, 2003.
\newblock In honour of Klaus-J\"{u}rgen Bathe.

\bibitem{GLOWINSKI1994283}
R.~Glowinski, T.-W. Pan, and J.~Periaux, ``A fictitious domain method for
  dirichlet problem and applications,'' {\em Comput. Methods Appl. Mech. Eng.},
  vol.~111, no.~3, pp.~283--303, 1994.

\bibitem{massing2014stabilized}
A.~Massing, M.~G. Larson, A.~Logg, and M.~E. Rognes, ``A stabilized nitsche
  fictitious domain method for the stokes problem,'' {\em Journal of Scientific
  Computing}, vol.~61, pp.~604--628, 2014.

\bibitem{SCHLOTTBOM2014}
M.~Burger, O.~Elvetun, and M.~Schlottbom, ``Analysis of the diffuse domain
  method for second order elliptic boundary value problems,'' {\em Foundations
  of Computational Mathematics}, 2014.

\bibitem{SCHLOTTBOM2016}
M.~Schlottbom, ``Error analysis of a diffuse interface method for elliptic
  problems with {D}irichlet boundary conditions,'' {\em Appl.~Numer.~Math.},
  vol.~109, pp.~109--122, 2016.

\bibitem{NEGRINI2021}
G.~Negrini, N.~Parolini, and M.~Verani, ``A diffuse interface box method for
  elliptic problems,'' {\em Applied Mathematics Letters}, vol.~120, p.~107314,
  2021.

\bibitem{SCHOTT2014233}
B.~Schott and W.~Wall, ``A new face-oriented stabilized {XFEM} approach for 2d
  and 3d incompressible {N}avier–{S}tokes equations,'' {\em Comput. Methods
  Appl. Mech. Eng.}, vol.~276, pp.~233--265, 2014.

\bibitem{BURMAN2010}
E.~Burman and P.~Hansbo, ``Fictitious domain finite element methods using cut
  elements: {I}. {A} stabilized {L}agrange multiplier method,'' {\em Comput.
  Methods Appl. Mech. Eng.}, vol.~199, no.~41-44, pp.~2680--2686, 2010.

\bibitem{BURMAN2012}
E.~Burman and P.~Hansbo, ``Fictitious domain finite element methods using cut
  elements: {II}. {A} stabilized {N}itsche method,'' {\em Appl.~Numer.~Math.},
  vol.~62, no.~4, pp.~328--341, 2012.

\bibitem{JASAKIBM2014}
H.~Jasak, D.~Rigler, and u.~Tukovi\'c, ``Design and implementation of immersed
  boundary method with discrete forcing approach for boundary conditions.''
  Proc. of 11th. World Congress on Computational Mechanics - WCCM XI, 2014.

\bibitem{VERSTEEG}
H.~K. Versteeg and W.~Malalasekera, {\em An introduction to computational fluid
  dynamics: the finite volume method}.
\newblock Harlow, 2007.

\bibitem{WELLER1998}
H.~G. Weller, G.~Tabor, H.~Jasak, and C.~Fureby, ``{A tensorial approach to
  computational continuum mechanics using object-oriented techniques},'' {\em
  Comput.~Phys.}, vol.~12, no.~6, pp.~620--631, 1998.

\bibitem{JASAK1996}
H.~Jasak, {\em Error Analysis and Estimation for Finite Volume Method with
  Applications to Fluid Flow}.
\newblock PhD thesis, Imperial College, University of London, 1996.

\bibitem{PATANKAR1980}
S.~V. Patankar, {\em Numerical Heat Transfer and Fluid Flow}.
\newblock Electro Skills Series, Hemisphere Publishing Corporation, 1980.

\bibitem{ELMAN2008}
H.~Elman, V.~Howle, J.~Shadid, R.~Shuttleworth, and R.~Tuminaro, ``A taxonomy
  and comparison of parallel block preconditioners for the incompressible
  {N}avier–{S}tokes equations,'' {\em J.~Comp.~Phys.}, vol.~227,
  pp.~1790--1808, 2008.

\bibitem{NEGRINIPHD}
G.~Negrini, {\em Non-conforming methods for the simulation of industrial
  polymer mixing processes}.
\newblock PhD thesis, Politecnico di Milano, 2023.

\bibitem{IKENO2007}
T.~Ikeno and T.~Kajishima, ``Finite-difference immersed boundary method
  consistent with wall conditions for incompressible turbulent flow
  simulations,'' {\em J.~Comp.~Phys.}, vol.~226, no.~2, pp.~1485--1508, Oct
  2007.

\bibitem{NEGRINI2023}
G.~Negrini, N.~Parolini, and M.~Verani, ``On the convergence of the
  {R}hie–{C}how stabilized {B}ox method for the {S}tokes problem,'' {\em
  Int.~J.~Numer.~Methods Fluids}, vol.~96, no.~8, pp.~1489--1516, 2024.

\bibitem{KIM2006}
N.~Kim, H.~Kim, and J.~Lee, ``Numerical analysis of internal flow and mixing
  performance in polymer extruder i: Single screw element,'' {\em Korea
  Aust.~Rheol.~J.}, vol.~18, pp.~143--151, 2006.

\bibitem{RAUWENDAAL2014}
C.~Rauwendaal, {\em Polymer Extrusion}.
\newblock Carl Hanser Verlag GmbH \& Company KG, 2014.

\end{thebibliography}

\end{document}